\documentclass[A4paper,12pt,epsf]{article}
\usepackage{amssymb}
\usepackage{amsfonts, eufrak}
\usepackage[dvips]{graphics}
\usepackage{amsmath}
\usepackage{amsthm}
\usepackage[latin1]{inputenc}
\usepackage[dvips]{graphics}

\newcommand{\R}{\mathbb {R}}
\newcommand{\Z}{\mathbb {Z}}
\newcommand{\N}{\mathbb {N}}
\newcommand{\E}{\mathbb {E}}
\newcommand{\I}{\mathbb {I}}
\newcommand{\prob}{\mathbb{P}}

\begin{document}

\centerline{\Large PARAMETER ESTIMATION IN }
\medskip

\centerline{\Large MANNEVILLE-POMEAU PROCESSES}
\thispagestyle{empty}

\vspace{0.8cm}

$${\bf Olbermann, B.P.}^{a}, {\bf Lopes, S.R.C.}^{b}
\footnote{ Corresponding author's E-mail: silvia.lopes@ufrgs.br},
\mbox{ and }   {\bf Lopes, A.O.}^b$$

\vspace{0.5cm}

\centerline{$^a$ Faculdade de Matem\'atica - PUCRS, Porto Alegre,
RS, Brazil}

\vspace{0.1cm}

\centerline{$^b$ Instituto de Matem\'atica - UFRGS, Porto Alegre,
RS, Brazil}

\vspace{0.5cm}

\centerline{July 10, 2007}

\vspace{0.5cm}

\abstract {\footnotesize{In this work we study a class of stochastic
processes $\{X_t\}_{t\in\N}$, where $X_t = (\varphi \circ
T_s^t)(X_0)$ is obtained from the iterations of the transformation
$T_s$, invariant for an ergodic probability $\mu_s$ on $[0,1]$ and a
continuous by part function $\varphi:[0,1] \to \R$. We consider here
$T_s:[0,1]\to [0,1]$ the Manneville-Pomeau transformation. The
autocorrelation function of the resulting process decays
hyperbolically (or polynomially) and we obtain efficient methods to
estimate the parameter $s$ from a finite time series. As a
consequence we also estimate the rate of convergence of the
autocorrelation decay of these processes. We compare different
estimation methods based on the periodogram function, on the
smoothed periodogram function, on the variance of the partial sum
and on the wavelet theory.}}

\vspace{0.5cm}


\vspace{0.3cm}


\newpage
\section{\normalsize Introduction}
\renewcommand{\theequation}{1.\arabic{equation}}
\setcounter{equation}{0}

The goal of this paper is to estimate the main parameter of some
processes obtained from iterations of Manneville-Pomeau maps.

We consider a class of stochastic processes $\{X_t\}_{t\in\N}$,
where $X_t = (\varphi \circ T_s^t)(X_0)$ is obtained from the
iterations of the transformation $T_s$, invariant for an ergodic
probability $\mu_s$ on $[0,1]$ and a continuous by part function
$\varphi:[0,1] \to \R$. The transformation $T_s:[0,1]\to [0,1]$ is
considered here as the Manneville-Pomeau map. We analyze the rate of
decay of the autocorrelation function for the resulting process. The
rate of convergence decays hyperbolically (or polynomially) not
exponentially. We obtain efficient methods to estimate the parameter
$s$ from a finite time series. As a consequence, we also estimate
the rate of convergence of the autocorrelation decay of these
processes.

Indeed, given $s$ the decay is known: Young (1999) has shown that
the autocorrelation decay of the Manneville-Pomeau processes has
order smaller than $n^{1-\frac1s}$, for $0<s<0.5$. Other models
which have similar properties to the Manneville-Pomeau map are the
linear by part approximation of the same map (see Fisher and Lopes,
2001 and Wang, 1989) and the Markov Chain with infinite symbols,
described in Lopes (1993).

Models of different phenomena in nature present autocorrelation
decay of the form $n^{-\beta}$, also called hyperbolic (or
polynomial): the use of the Markov Chain model seems to be
appropriated for the analysis of DNA sequences (see Peng et al.,
1992 and 1996 and Guharay et al., 2000); cardiac rhythm fluctuations
(see Absil et al., 1999 and Peng et al., 1996); turbulence (see
Schuster, 1984) and economy (see Mandelbrot, 1997; Lopes et al.,
2004 and Lopes, 2007). In most of the cases the exact rate of
convergence of the autocorrelation function decay is a relevant
information in the model. Here we are interested in to compare
different methods for estimating such $\beta$ in the case of the
Manneville-Pomeau processes.

When $0.5<s<1.0$ we have the {\it long range dependence regime\/}.
{\it Fractionally integrated autoregressive moving average\/}
(ARFIMA) models also present such behavior (see Beran, 1994; Geweke
and Porter-Hudak, 1983; Reisen and Lopes, 1999 and Lopes et al.,
2004). The corresponding parameter for the ARFIMA model is $d=1-
\frac{1}{2s}$. The ARFIMA process has an explicit formula for the
spectral density function $f_X(\cdot)$ (see Reisen et al., 2001;
Lopes et al., 2002 and Olbermann et al., 2006) but this is not the
case for the processes considered here.

When $0< s <0.5$ we have the {\it not so long dependence regime\/}.
The so-called {\it intermediate dependence regime} happens when
$s\in (\frac13,\frac12)$.

Recently several interesting papers appear describing the statistics
of time series obtained from dynamical systems: Chazottes (1998),
Chazottes (2005), Collet et al. (1995), Collet et al. (2004) and
Collet (2005). We also refer the reader to the last sections of the
book by Collet and Eckmann (2006).

Here we analyze and compare several estimation procedures based on
the periodogram function, on the smoothed periodogram function, on
the variance of the partial sum and on wavelet theory.

The paper is organized as follows. In Section $2$ we define the
Manneville-Pomeau maps giving some definitions, basic properties and
results. Section $3$ presents the Manneville-Pomeau processes that
will be the setting of the estimation procedures we choose in this
work. In Section $4$ we consider the estimation procedures for the
long dependence case while in Section $6$ we present the Monte Carlo
simulation study for this regime. In Section $5$ we consider the
estimation procedures for the not so long dependence case whlie in
Section $7$ we present the Monte Carlo simulation study for this
other regime. Section $8$ contains a summary of the paper. In
Appendix $A$ we consider some general properties of the Fourier
series which are necessary in the paper. Appendix $B$ contains the
theoretical reasoning for some of the estimation procedures proposed
in Section $4$ of the paper.

\vspace{0.3cm}

\section{\normalsize Manneville-Pomeau Maps}
\renewcommand{\theequation}{2.\arabic{equation}}
\setcounter{equation}{0}

In this section we present the Manneville-Pomeau maps, some
definitions, basic properties and results.

We first define the Manneville-Pomeau transformation
and we give some of its properties.

\vspace{0.3cm}

\noindent
{\bf Definition 2.1}: Let $T_s: [0,1]\to[0,1]$ be the {\it
Manneville-Pomeau map\/} given by

\begin{eqnarray}\label{MPmaps}
T_s(x) = x+ x^{1+s} \pmod{1} = \left\{
\begin{array}{lcl}
x+x^{1+s}, & \mbox{ if } & x+x^{1+s}\leq 1\\
x+x^{1+s}-1, & \mbox{ if } & x+x^{1+s} > 1,
\end{array}\right.
\end{eqnarray}

\noindent
where $s$ is a positive constant.

As usual, we shall use the following notation
$$
T_s^t\equiv\underbrace{T_s\circ \cdots \circ T_s}_{t-times}.
$$

\vspace{0.3cm}

The map $T_s$ (see Figure $2.1$ (a)), given by the expression
(\ref{MPmaps}) has the following properties:

\begin{itemize}

\item $T_s$ is a piecewise monotone function with two
full branches, that is, there exists $p\in \N-\{0\}$ such that
$T_s|_{(0,p)}$ and $T_s|_{(p,1)}$ are strictly monotone, continuous
and $T_s((0,p))=(0,1)=T_s((p,1))$, where $p+p^{1+s}=1$.

\item The branches $T_s|_{(0,p)}$ and $T_s|_{(p,1)}$
are $C^2$.

\item $T_s'(x)>1$, for all $x>0$, and
$T_s'(x)\geq \lambda>1$, for $x\in (p,1)$.

\item $T_s$ has a unique indifferent fixed point $0$.
Therefore, $T_s(0)=0$ and $|T_s'(0)|=1$.

\item There exists an invariant absolutely continuous
ergodic measure
$\mu_s$ for the Manneville-Pomeau transformation $T_s$.
Thaler (1980), using the properties of the
Manneville-Pomeau map, shows that
$d\, \mu_s(x)\equiv h_s(x)\, dx$,
where $
h_s(x) \approx x^{-s}, {\mbox { for }}
x\in (0,1)$, close to 0.
\noindent

When $s\geq 1$, the measure $\mu_s$ has infinite mass and it is not
a probability.

When $0<s<1$, the probability $\mu_s$ is {\it mixing} for $T_s:
[0,1] \to [0,1]$ (see Young, 1999; and Fisher and Lopes, 2001).

\end{itemize}

\vspace{0.5cm}

Given a continuous by part function $\varphi:[0,1]\to \R$, one can
consider the random variables $X_t=(\varphi \circ T_s^t )(X_0)$, for
$t\in \N$, where $X_0$ is distributed according to the probability
$\mu_s$. The stationary stochastic process $\{X_t\}_{t\in \N}$ is
called the {\it Manneville-Pomeau process\/}. We will consider here
$\varphi$ as an indicator function of an interval in $[0,1]$. In
this case, the time series obtained from the process $\{X_t\}_{t\in
\N}$ will be a binary time series of 0's and 1's only.

It is known that the autocorrelation decay of the Manneville-Pomeau
processes, given by the expression (\ref{MPprocess}), has order
smaller than $n^{1-\frac1s}$, for $0<s<0.5$ (see Young, 1999). In
Fisher and Lopes (2001) it is shown, for the linear by part model
given by Definition 2.2 below, that these bounds are exact (for the
corresponding values).

We refer the reader to Maes et al. (1999) for more details on the
dynamics of the system given by (\ref{MPmaps}).

Other models which have similar properties to the Manneville-Pomeau
map are the linear by part approximation of the same map (see
Definition 2.2 below and Fisher and Lopes, 2001 and Wang, 1989) and
the Markov Chain with infinite symbols (see Definition 2.3 below)
described in Lopes (1993). The use of the Markov Chain model
$\{Y_t\}_{t\in \N}$, defined below, seems to be appropriated for the
analysis of DNA sequences (see Peng et al., 1992 and 1996). The same
estimation methods, proposed for the Manneville-Pomeau processes in
Section $4$, can be also applied to these other models.

\vspace{0.3cm}

\noindent {\bf Definition 2.2}:  Let $\zeta(\gamma)=\sum_{n\geq 1}
n^{-\gamma}$ be  the {\it Riemann zeta} function. Consider the
partition in intervals of $[0,1]$ given by
$$
M_0=\left (1- \frac{1}{\zeta(\gamma)},1\right ){\mbox { and}}\,
M_k = \left (1- \frac{1}{\zeta(\gamma)}\sum_{n=1}^{k-1}
n^{-\gamma},1- \frac{1}{\zeta(\gamma)}\sum_{n=1}^{k}
n^{-\gamma}\right ),
$$
\noindent for $k\geq 1$. For $\gamma >2$, we define the following
linear by part transformation $T_{\gamma}:[0,1]\to [0,1]$ such that
over the interval $M_k$, for $k\geq 1$, $T_{\gamma}$ has slope
$((k+1)k^{-1})^{\gamma}$ and over the interval $M_0$ it has slope
${\zeta(\gamma)}$. We assume that the branches
$$
T_{\gamma}|_{\left( 0,1- \frac{1}{\zeta(\gamma)}\right)} {\mbox {
and }} T_{\gamma}|_{M_0}
$$
\noindent are continuous; under these assumptions the transformation
$T_{\gamma}$ is uniquely defined (see Figure $2.1$ (b)). The
transformation $T_{\gamma}$ is called {\it the linear by part
approximation of the Manneville-Pomeau map}.

\vspace{0.3cm}

In the same way as before, giving a continuous by part function
defined by $\varphi:[0,1]\to \R$, one can consider the random
variables $X_t=(\varphi \circ T_{\gamma}^t )(X_0)$, for $t\in \N$,
where $X_0$ is distributed according to a certain probability
$\mu_{\gamma}$, invariant for $T_{\gamma}$. The probability
$\mu_{\gamma}$ is absolutely continuous with respect to the Lebesgue
measure. We call $\{X_t\}_{t\in \N}$ the {\it linear by part
approximation of the Manneville-Pomeau process\/}.

Each value of  $s$ for the Manneville-Pomeau map corresponds to a
value $\gamma= 1+\frac1s$ with the same behavior with respect to the
autocorrelation decay.

\vspace{0.3cm}


\vspace{0.3cm}

The Manneville-Pomeau map has the advantage of been more suitable
than the linear by part model for computer implementation when one
is interested on Monte Carlo simulations. For this reason, in the
simulation sections we will concentrate our analysis in such model.

Below we define a Markov process with state $\N$ based on a certain
transition probability matrix ${\bf P}$. The time evolution of such
process will also have similarities with the iteration of
Manneville-Pomeau map.

\vspace{0.3cm}

\noindent {\bf Definition 2.3}: Let {\bf P} be a Markov chain with
infinite transition probability matrix {\bf
P}=($\prob(i,j))_{i,j\in\N}$ (see page 153 in Lopes, 1993; Wang,
1989 and Feller, 1949) with transition probabilities given by
$$
\prob(n,n-1) =1,\, {\mbox{ for all }} n\in\N-\{0\},
$$
$$
\prob(n,j)= 0, \, {\mbox{ for }} j\neq n-1,
$$
\noindent
and
$$
\prob(0,n) = \frac{(n+1)^{-\gamma}}{\zeta (\gamma)},
$$
\noindent where $\zeta(\gamma)$ is the Riemann zeta function and
$\gamma > 2$. There exists an explicit formula for the eigenvector
$\pi_0$ associated to the eigenvalue 1 (see page 154 in Lopes,
1993).

\vspace{0.3cm}

Let $\{Z_t\}_{t\in \N}$ be the stationary stochastic Markov process
obtained from the transition matrix {\bf P} above and from the
initial stationary distribution $\pi_0$. Let $\I_0$ be the indicator
function of the set $A=\{0\}$ on $\N$. Let now $\{Y_t\}_{t\in \N}$
be the process $1 - \I_0(Z_t)$. In this way, we identify paths
$\omega \in \N^{\N}$ with paths $\tilde \omega \in \{0,1\}^{\N}$.
Then, $\{Y_t\}_{t\in \N}$ is a stochastic process with random
variables assuming only the values $0$ and $1$. For the process
$\{Y_t\}_{t\in \N}$ consider the probability induced by the process
$\{Z_t\}_{t\in \N}$ by means of the identification of the paths.

To clarify the ideas in the above Definition $2.3$, the following
example shows the identification paths in $\N^{\N}$ to paths in
$\{0,1\}^{\N}$.

\vspace{0.5cm}

\noindent \bf Example 2.1: \rm Let $\{Z_t\}_{t\in\N}$ be the process
where a sample path $w \in {\N}^{\N}$, for instance,
$w=\{0765432109876543210543210\cdots\}$, is associated to another
sample path of the process $\{Y_t\}_{t\in\N}$. The corresponding
sample path for the process $\{Y_t\}_{t\in\N}$ is given by
$$
\tilde w = \{0{\underbrace{1111111}_7}0 {\underbrace{111111111}_9}
0{\underbrace{11111}_5}0\cdots\}.
$$
\noindent Hence, we applied the change of coordinates $Z_t \to Y_t$
associating sequences of natural numbers to blocks of $1$
intercalated by $0$, in such a way that the structure of the process
is kept the same.

\vspace{0.5cm}

We say that two different stochastic processes are {\it equivalent}
when there is a bijective change of coordinates acting in the set of
paths transferring the probability of one process into the other.

The process $\{Z_t\}_{t\in \N}$ is, by definition, {\it equivalent}
to the process $\{Y_t\}_{t\in \N}$ by the above change of
coordinates. One can also show that $Y_t$ is also {\it equivalent}
to $X_t=(\varphi\circ T^t_{\gamma})(X_0)$ (see Section $4$ in Lopes,
1993 with $\varphi\equiv {\I}_{M_0}$).

It is also known that the central limit theorem (converging to a
Gaussian distribution) is true for the Manneville-Pomeau stochastic
process $\{X_t\}_{t\in \N}$, described in Section $3$, when $0<
s<0.5$ due to the rate of convergence of the autocorrelation decay
(see Young, 1999; Lopes, 1993 and pages 1099-1100 in Fisher and
Lopes, 2001).

When $0.5<s<1.0$ it was conjectured that for the Manneville-Pomeau
stochastic process $\{X_t\}_{t\in \N}$ the central limit theorem is
true, but it converges to a stable law with parameter
$\alpha=s^{-1}$. This was proved by Gou\"ezel (2003).  From Feller
(1949) it is known for the corresponding parameter of the Markov
Chain model described above (or for the equivalent process
$X_t=(\varphi\circ T^t_{\gamma})(X_0)$ with $\varphi\equiv
{\I}_{M_0}$ (see Wang, 1989 or Section $4$ in Lopes, 1993, for more
details)).

For the estimation in the {\it long range dependence case\/}, one
has to consider larger sample sizes for the time series. In this
situation, in general, the computation effort for obtaining good
results is very high. This is something that one can not avoid due
to the small rate of convergence of the decay. The mixing rate is
not so good as it happens, for instance, when one considers models
with exponential autocorrelation decay. We present here several
quite efficient methods to obtain reasonable results. One method is
by using the periodogram function described in Sections $4$ and $6$.
The method based on wavelet works fine in several cases and
surprisingly can also be applied to estimate $s$ when $s\geq 1.0$
(see Sections $4$ and $6$).

In a forthcoming paper by A.S Pinheiro and S.R.C. Lopes will be
presented a bias correction for the wavelet estimation in {\it
long\/} and {\it not so long dependence cases}.

\vspace{0.3cm}

\section{\normalsize Manneville-Pomeau Process and Some of its Properties}
\renewcommand{\theequation}{3.\arabic{equation}}
\setcounter{equation}{0}

In this section we define the Manneville-Pomeau stochastic processes
and present some of their properties.

\vspace{0.3cm}

Let  $\varphi:[0,1]\to \R$ be a $\mu_s$-integrable function and
$T_s(\cdot)$ the Manneville-Pomeau transformation given by the
expression (\ref{MPmaps}). The {\it  Manneville-Pomeau stochastic
process} $\{X_t\}_{t\in\N}$ is given by
\begin{equation}\label{MPprocess}
X_t =(\varphi \circ T_s^t )(X_0)= \varphi ( T_s^t (X_0))=
\varphi ( T_s(X_{t-1}))= (\varphi \circ T_s)(X_{t-1}),
\end{equation}
for all  $t\in \N$, where $X_0$ is distributed  according to the
measure $\mu_s$. In other words, the Manneville-Pomeau process
$\{X_t\}_{t\in\N}$ is obtained applying  $\varphi$ to the
iterations of $T_s$, that is, $X_t=\varphi \circ T_s^t$, for $s$
fixed and  $t\in\N$.

\vspace{0.5cm}

We shall consider here only the case where $\varphi$ is the
indicator function ${\I}_{A}$ of an interval $A$ contained in
$[0,1]$ or else $\varphi= {\I}_{A} - \mu_s (A)$. Our simulations,
shown in Sections $5$ and $7$, will be done for the case where $A=
[0.1,0.9]$.

We shall denote by $\gamma_X(\cdot)$ the autocovariance function
for the process $\{X_t\}_{t\in \N}$, that is,
\begin{equation}\label{autocov}
\gamma_X (h) \equiv {\E}_\mu (X_h X_0) - [{\E}_\mu(X_0)]^2= \int
\varphi( T^h (x)) \varphi(x)  d \mu_s (x) - [\int \varphi (x) d
\mu_s (x) ]^2,
\end{equation}
\noindent
$\mbox{ for } h\in \N$.

We denote by $\rho_X(\cdot)$ the autocorrelation function
of the process $\{X_t\}_{t\in \N}$, that is,
\begin{eqnarray*}
\rho_X(h)=\frac{\gamma_X(h)}{\gamma_X(0)},
\, \mbox{ for all } h\in \N,
\end{eqnarray*}
where $\gamma_X(0)\equiv {\E}_\mu (X_0^2) - [{\E}_\mu(X_0)]^2=
Var_\mu(X_0)$ is the variance of the process.

\vspace{0.3cm}

The spectral density function of the process
$\{X_{t}\}_{t\in \N}$ is given by
\begin{equation}\label{spectral}
f_X(\omega)=\frac{1}{2\pi}[\gamma_{X}(0)+2\sum_{h=1}^{\infty}
\gamma_{X}(h)\cos(\omega h)],\mbox{ for } \omega\in[-\pi,\pi].
\end{equation}

Now we shall define the {\it periodogram function\/}
associated to a time series $T_s^t (x_0)$, for
$1 \leq t \leq N$, obtained from a $x_0$ chosen with
probability one according to
the measure $\mu_s$. The periodogram function is given by
\begin{equation}\label{periodogram}
I(\omega_h) = f_N(\omega_h) \overline{f_N(\omega_h)},
\end{equation}
where
\begin{eqnarray*}
f_N (\omega) = \frac{1}{2\pi \sqrt{N}} \sum^N_{t=1}
\varphi(T_s^t(x_0)) e^{-i\omega t}, \ \ \omega\in (0,2\pi],
\end{eqnarray*}
\noindent
with $\overline{f_N(\cdot)}$ indicating the complex
conjugate of $f_N(\cdot)$ and
\begin{equation}\label{omegah}
\omega_h = \frac{2\pi h}{N}, \mbox{ for } h=0,1, \cdots, N,
\end{equation}
the $h$-th discrete Fourier frequency (see Brockwell and Davis,
1991).

Note that the periodogram function depends on $x_0$ and $N$ (large).
One can obtain a good approximation of the spectral density function
$f_X(\cdot)$ by the periodogram function (see Lopes and Lopes, 2002
for a mathematical proof that can be applied to the case we analyze
here when $0<s<0.5$).

The {\it periodogram function\/} is an unbiased estimator for the
spectral density function $f_X(\cdot)$, even though it is not
consistent (see Brockwell and Davis, 1991).

Another procedure for estimating the parameters which produces good
results is by using the wavelet theory. This type of analysis can be
also used in the regime $s>1$ where the spectral density function,
defined in the expression (\ref{periodogram}), does not exist since
the random process is not associated to a probability.

\vspace{0.5cm}

We shall use the following notation:

\vspace{0.2cm}

\begin{itemize}

\item{} If, for the sequence $\{a_n\}_{n\in\N}$, there
exists $u\in \R$ and, for any $\delta >0$, there
exist positive constants $c_1$ and $c_2$
such that, for all $n\in \N$,
$$
c_1\, n^{u-\delta}\leq \left | a_n \right |
\leq c_2\, n^{-u+\delta},
$$
then we denote $a_n\approx n^{-u}$. We also
say that $a_n$ is {\it of order\/} $n^{-u}$,
for $n\to \infty$.

\item{} If, for the real function  $g(\cdot)$, there
exist $b\in \R$ and $\epsilon >0$ such that, for any $\delta >0$,
there exist positive constants
$d_1$ and $d_2$ such that, for all $x\in (0,\epsilon)$,
$$
d_1\, x^{b+\delta}\leq \left |g(x) \right |
\leq d_2\, x^{b-\delta},
$$
then, we denote $g(x)\approx x^{b}$. We also say
that $g$ is {\it of order\/} $x^b$ around $0$.

\end{itemize}

\vspace{0.3cm}

If there exist $c_1,c_2>0$ such that
$$
c_1\, n^{-u}\leq \left | a_n \right |
\leq c_2\, n^{-u},
$$
then, of course, $a_n\approx n^{-u}$.

If there exist $d_1,d_2>0$ such that
$$
d_1\, x^b \leq \left |g(x) \right |
\leq d_2\, x^b,
$$
then, of course, $g(x)\approx x^{b}$. We need
however, this more general definition because of
Theorem A.4 in the Appendix A of the present
work.

\vspace{0.3cm}

\noindent
{\bf Definition 3.1}: Let $\{X_{t}\}_{t\in\N}$ be a
 stochastic stationary process with autocovariance
function $\gamma_X(\cdot)$ given by the expression (\ref{autocov}).
If there exists $u \in (0,1)$ such that
\begin{equation}
 \gamma_{X} (h)\approx{h^{- u}},
\end{equation}
\noindent
then we say that $\{X_{t}\}_{t\in \N}$ is a
{\it stochastic process with long dependence}.

\vspace{0.3cm}

\noindent {\bf Definition 3.2}: Let $\{X_{t}\}_{t\in \N}$ be a
stochastic stationary process with autocovariance function
$\gamma_X(\cdot)$ given by the expression (\ref{autocov}). If there
exists $u>1$ such that
\begin{equation}
 \gamma_{X} (h)\approx{h^{- u}},
\end{equation}
\noindent
then we say that $\{X_{t}\}_{t\in \N}$ is a
{\it stochastic process with not so long dependence\/}.

\vspace{0.3cm}

For the Manneville-Pomeau process it is known that
\begin{equation}
\gamma_X (h)\approx h^{1-\frac1s},
\end{equation}
(see Young, 1999 for the upper bound and Fisher and Lopes, 2001 for
the lower bound).

\vspace{0.3cm}

When $0.5<s<1$ the Manneville-Pomeau process, given by the
expression (\ref{MPprocess}), has the {\it long dependence\/}
property and when $0<s<0.5$ it has the {\it not so long
dependence\/} property. We shall consider here different methods for
estimating the value $s$ in both cases.

\vspace{0.2cm}

In the {\it long dependence regime\/} there exists a relationship
between the velocity of the autocorrelation function decay to zero
and the regularity of the function $f_X(\cdot)$. This property
follows just from a careful analysis of Fourier series. We refer the
reader to Chapter X, Section $3$ in Bary (1964), pages 1086-1090 in
Fisher and Lopes (2001) and also the Appendix A of the present work
for a careful description of this relationship. This follows
basically from the fact that if $f_X(\lambda)\approx \lambda^{-b}$,
with $b>0$, then $\gamma_X (h) \approx h^{b-1}$. In the case when
the coefficients $\gamma_{X} (h)$ are monotone decreasing in $h$,
then $f_X(\lambda)\approx \lambda^{-b}$, if $\gamma_X(h)\approx
h^{b-1}$, for $b>0$. Fisher and Lopes (2001) show that the
autocovariance function $\gamma_X (h)$ is a monotone function for
the linear by part approximation of the Manneville-Pomeau map in the
case of a certain $\varphi$. These authors also show that $\gamma_X
(h)\approx{h^{\gamma - 3}}$, when $2<\gamma<3$ (see page 1090).

In the case of Manneville-Pomeau maps with {\it long dependence\/},
from the exact asymptotic given by the expression (\ref{spectral}),
one can obtain (by analogy with the linear by part model) the rate
of convergence of the autocorrelation decay to zero from the
asymptotic of $f_X(\lambda)$ to infinity when $\lambda \to 0$ and
vice versa. It follows from the above considerations and from
(\ref{spectral}) that $f_X(\lambda)\approx \lambda^{\frac1s -2}$.

The phenomena $f_X(\omega)\approx \omega^{-b}$ is known as
$\frac{1}{f}$-noise property (in this case, $\frac{1}{f^b}$-noise
would be a more appropriate terminology), where $f$ stands for a
frequency (here denoted by $\omega$).

\vspace{0.3cm}

\noindent
{\bf Definition 3.3}: The continuous function $g: (-\pi,
\pi)\to {\R}$ is said to be {\it H\"older of order\/} $a$, $0<a<1$,
if there exists a positive constant $K$ such that
$$
|g(x) - g(y)| \leq K |x-y|^a,
$$
for any $x,y \in (-\pi, \pi)$. We also call
$a$ the {\it exponent of\/} $g$.

\vspace{0.3cm}

\noindent
{\bf Definition 3.4}: The continuous function $g: (-\pi,
\pi)\to {\R}$ is said to be exactly $a$-{\it H\"older} in the point
$x_0$, for $0<a<1$, if for any $\delta>0$, there exist positive
constants $c_1$ and $c_2$ such that
$$
c_1\, |x-y|^{a+\delta} \leq |g(x) - g(y)|
\leq c_2 \, |x-y|^{a-\delta},
$$
for any $y \in (-\pi, \pi)$. We also call
$a$ the {\it exact exponent of\/} $g$ at
$x_0$.

We will apply this definition for the case $x_0=0$.

\vspace{0.3cm}

When one considers the Manneville-Pomeau maps with {\it not so long
dependence\/}, one can say more about the regularity of $f_X(\cdot)$
(see Chapter II, Section $3$ and Chapter X, Section $9$ in Bary,
1964 and Appendix A of this present work): it is exactly
$\beta$-H\"older continuous function with exponent $\beta=\frac1s
-2$. We are using here the notation: a $\beta$-H\"older function,
with $\beta =n+\alpha$, $0<\alpha<1$, is a function such that it is
$n$ times differentiable and the $n$-th derivative is
$\alpha$-H\"older.

The periodogram function $I(\cdot)$ is a useful way to obtain an
approximation of $f_X(\cdot)$ (see Lopes and Lopes, 2002). One can
obtain an estimation of $s$ from the above considerations and from
the periodogram function as we will explain in the next section.

\vspace{0.3cm}

\section{\normalsize Estimation in the ``Long Dependence" Case}
\renewcommand{\theequation}{4.\arabic{equation}}
\setcounter{equation}{0}

The main goal of this section is to estimate the transformation
$T_s$, or equivalently, to estimate the parameter $s$, when $0.5 <s<
1$. For this purpose, we consider a finite time series
$\{X_t\}_{t=0}^{N-1}$ obtained from the process $\{X_t\}_{t\in\N}$
given by (\ref{MPprocess}).

By Monte Carlo simulation, that is given in Section $5$, we compare
some methods for estimating $s$ with the one presented in Schuster
(1984). We are interested in the performance of this method when
compared to the others.

The process $\{X_t\}_{t\in\N}$, defined by the expression
(\ref{MPprocess}), is considered here to be
\begin{equation}\label{MPprocess1}
X_t={\I}_{A} \circ T_s^t={\I}_{(0.1,0.9)} \circ T_s^t,
\end{equation}
which is stationary and ergodic (see Lopes and Lopes, 1998).

\vspace{0.3cm}

For the {\it long dependence case\/} one can express the graph of
$f_X(\cdot)$ (or of the periodogram function $I(\cdot)$) in the
logarithm scale, and this exhibits linear behavior. By ordinary
least squares estimation one can obtain an estimate of the value
$s$.

We now explain more carefully this very useful method for the {\it
long dependence case\/}: suppose there exists  $c$ such that
$f_X(\omega) \approx \omega^c$, for $\omega$ close to zero. Then,
for $\omega$ close to zero
$$
\frac{\ln (f_X(\omega))}{\ln (\omega)} \approx c.
$$

\noindent
From the estimated value of $c$ we estimate $s$
since $c= \frac1s -2$. An estimate of $c$ can be obtained
via the periodogram by
$$
\frac{\ln (I(\omega))}{\ln (\omega)} \approx \hat c,
$$
\noindent with $\omega$ chosen very close to $0$.

We shall now consider six different methods for estimating the
parameter $s$: the least squares method proposed in section $4.3$ of
Schuster (1984); the least squares method proposed here using the
smoothed periodogram function when the Parzen or the ``cosine bell"
lag window are used to consistently estimate the spectral density
function; the one based on the variance of the sample partial sums
of the process; the one based on the logarithm of the variance of
the sample mean of the process and the one based on the wavelet
theory. These methods are described in this section and in Section
$5$ we present a Monte Carlo simulation study comparing them.

\vspace{0.5cm}

\noindent {\bf {\it Perio} Estimator }

\vspace{.5cm}

This method is based on the periodogram function of a time series
$\{X_t\}_{t=1}^N$ and it is largely used by the physicists (see
Schuster, 1984).

The estimator of $s$ is obtained from the least squares method based
in a linear regression of $y_1,y_2,\cdots,y_{g(N)}$ on
$x_1,x_2,\cdots,x_{g(N)}$, where $y_j=\ln (I(\lambda_j))$,  $x_j=\ln
(j)$ and  $g(N)=N^{0.5}$. The $I(\cdot)$ is the periodogram function
given by the expression (\ref{periodogram}) and $\lambda_j$ is the
$j$-th Fourier frequency given by (\ref{omegah}). Let $c$ be the
slope coefficient of the linear regression in the logarithm scale.
The coefficient $c$ allows the estimation of $s$ through the
equality
$$
s=\frac {1}{c+2},
$$
since, for $s\in (0,1)$ we know that

$$
f_X(\omega) \approx \omega^{\frac{1}{s}-2}, {\mbox { for }} \omega
{\mbox { close to the zero frequency}}.
$$
\noindent
Therefore,
\begin{equation}\label{perio}
{\hat c}=\frac {1}{\hat s}-2
\Leftrightarrow {\hat s}=\frac{1}{\hat c+2}.
\end{equation}
\noindent We shall denote the estimator in (\ref{perio}) by $Perio$.

\vspace{0.5cm}

\noindent {\bf {\it Parzen} Estimator }

\vspace{0.5cm}

This method is also a regression estimator for the parameter $s$ and
is obtained by replacing the periodogram function $I(\cdot)$ in the
{\it Perio} method by its smoothed version with the Parzen lag
window (see Brockwell and Davis, 1991). It is known that the use of
a spectral lag window consistently estimates the spectral density
function (see Brockwell and Davis, 1991). This estimator has the
same expression as in (\ref{perio}), but now
$y_j=\ln({f_{sm}(\omega_j))}$, where $f_{sm}(\cdot)$ is the smoothed
periodogram function. The value of $g(N)$ is chosen as in the {\it
Perio} method. The truncation point in the Parzen lag window is
considered to be $m=N^{0.9}$.

\vspace{0.5cm}

\noindent {\bf {\it Cos} Estimator }

\vspace{0.5cm}

This method is similar to the {\it Parzen} estimator, where now one
uses the ``cosine bell" spectral lag window (see Brockwell and
Davis, 1991). Its expression is given by (\ref{perio}), where now
the smoothed periodogram function $f_{sm}(\cdot)$ is obtained from
the ``cosine bell" lag window. Again, by a linear regression we
obtain the estimator of $s$. In this method we considered different
limits for $g(N)=N^{\alpha_i}$: we used $\alpha_1= 0.5$ and
$\alpha_2=0.7$ and we denote this estimator by $Cos(i)$, $i=1,2$.

\vspace{0.3cm}

\noindent {\bf Remark 4.1:} The methods $Perio$, $Parzen$ and $Cos$,
defined above, are similar to those proposed by Lopes et al. (2004)
and Reisen et al. (2001) to estimate the differencing parameter in
ARFIMA models. They are also similar to the estimators proposed by
Lopes (2007) for the differencing $d$ or the seasonal differencing
$D$ parameters in seasonal fractionally integrated
ARIMA$(p,d,q)\times(P,D,Q)_s$ process with period $s$. Again, we
observe that there is no exact expression for the spectral density
function $f_X(\cdot)$ in the case of the Manneville-Pomeau
processes.

\vspace{0.5cm}

\noindent {\bf {\it Varmp} Estimator }

\vspace{0.5cm}

This method, denoted by $Varmp$, is different from the other
previous three. To explain this method, we consider a time series of
sample size $N$ from the process (\ref{MPprocess1}) and let $M_N$ be
the random variable given by
\begin{equation}\label{partialsum}
M_N={\mbox{ total number of  1's in the time series }}
\{X_t\}_{t=0}^{N-1} =
\sum_{i=0}^{N-1}X_i=
S_N.
\end{equation}
\noindent One can show (see Lopes, 1993; Olbermann, 2002 or Wang,
1989) that
\begin{equation}\label{varpartialsum}
Var(M_N)\approx N^{4 -\gamma}= N^{3-\frac1s}.
\end{equation}
\noindent
We present a proof of this fact in a quite large
generality in Appendix B.

The property (\ref{partialsum}) allows one to obtain another
estimator for the parameter $s$. In fact, if one applies the
logarithm to that expression one gets
$$
Varmp =\frac{1}{3-\frac{\ln(Var(M_N))}{\ln (N)}}=
\hat s.
$$

\vspace{0.3cm}

\noindent {\bf Remark 4.2:} As in the ARFIMA process (see Beran,
1994 and Olbermann, 2002) we observe that this estimator is also
very much biased to estimate $s$ in the Manneville-Pomeau processes.

\vspace{0.5cm}

\noindent {\bf {\it Vpmp} Estimator }

\vspace{0.5cm}

This method is also based on the variance of the random variables
$M_N$. It is proposed by Beran (1994) under the name of {\it
variance plot\/}. It is obtained from the order of the variance of
${\bar X_N}=\frac{S_N}{N}$ given by
\begin{equation}\label{varmean}
\mbox{Var}({\bar X_N})\approx O(N^{2d -1}),
\end{equation}
\noindent
where $d$ is the differencing parameter
in ARFIMA models.

For the Manneville-Pomeau processes we only need to consider the
expression (\ref{varmean}), the relationship between the random
variables $M_N$ and $S_N$, given by (\ref{varpartialsum}) and the
relationship between the parameters $s$ and $d$, given by
$d=1-\frac{1}{2s}$. We shall denote this estimator by $Vpmp$.

\vspace{0.5cm}

\noindent {\bf {\it Wmp} Estimator }

\vspace{0.5cm}

This method is based on the wavelet estimator
proposed by Jensen (1999) to estimate the differencing
parameter $d$ in ARFIMA models. To consider this
a a method to estimate the parameter
$s$ in Manneville-Pomeau processes we must
consider the relationship between the parameters
$s$ and $d$, given by $d=1-\frac{1}{2s}$ and
the estimator proposed here, denoted by $Wmp$.

We refer the reader to Percival and Walden (1993) and Lopes and
Pinheiro (2007) for the use of wavelets in several different
problems in statistics.

A {\it wavelet\/} is any continuous function $\psi(t)$ that decays
fast to zero when $|t|\to \infty$ and oscillates in such a way that
$\int_{-\infty}^ \infty \psi(t)\, dt =0$. The idea is to use diadic
translations and dilations of the function $\psi(\cdot)$ such that
they generate the whole ${\cal L}^2(\R)$. From this, the wavelets
considered are of the form
$$
\psi_{j,k}(t)=2^{\frac{j}{2}}\, \psi(2^j\, t - k),
\, \mbox{ for } \, j, k \in \Z,
$$
\noindent which constitute an orthonormal basis of ${\cal L}^2(\R)$
(see Percival and Walden, 1993).

Here we consider only the wavelet bases Haar and Mexican
hat, since these bases have easy analytic expressions
given by
\begin{eqnarray*}
\psi_{j,k}(t)=\left\{\begin{array}{lll}2^{\frac{j}{2}},
& \mbox{ if } 2^{-j}\, k\leq t<2^{-j}\, (k+\frac12)\\
-2^{-j}, & \mbox{ if } 2^{-j}\, (k+\frac12)\leq t <
2^{-j}\, (k + 1) \\
0, & \mbox{ otherwise }
\end{array} \right.
\end{eqnarray*}
and
$$
\psi_{j,k}(t)=2^{\frac{j}{2}}[1-(2^j\, t-k)^2]
\exp[-(2^j\, t-k)^2/2],
$$
for $j=0,1,\cdots,m-1$ and $k=0,1,\cdots, 2^j-1$,
where $m\in N$ is such that $N=2^m$.

Given a time series of the sample size $N$ from the stochastic
process (\ref{MPprocess1}) we define the {\it wavelet
coefficients\/} as the finite wavelet transform for this time series
given by
$$
\omega_{j,k}= 2^{\frac{j}{2}}\sum_{t=0}^{N-1}
X_t\, \psi(2^j\, t-k),
$$
for $j=0,1,\cdots,m-1$ and $k=0,1,\cdots, 2^j-1$,
where $m\in N$ is such that $N=2^m$.

To obtain the estimator proposed by Jensen (1999)
we define the {\it variance of the wavelet
coefficients\/} as
$$
R(j)=\E[(\omega_{j,k})^2], \, \mbox{ for all }
\, j=0,1,\cdots, m-1.
$$

Considering the relationship between
$s$ and $d$ given by
$d=1-\frac{1}{2s}$,
the estimator based on the wavelets
is given by
$$
Wmp =\frac{\sum^{m-1}_{j=4} x_j^2}{2\, \left( \sum_{j=4}^{m-1} x_j^2
- \sum^{m-1}_{j=4} x_j\ln (\hat R(j)) \right)},
$$
\noindent
where $x_j$ is given by
$$
x_j = \ln (2^{-2j}) - \frac{1}{m-4}\sum_{j=4}^{m-1} \ln (2^{-2j}),
$$
\noindent
and $\hat R(j)$ is the {\it sample variance of the
wavelet coefficients\/} defined by
$$
\hat R(j)\equiv \frac{1}{2^j}\sum_{k=0}^{2^j -1}
(\omega_{j,k})^2, \ \mbox{ for all }
\, j=4,5,\cdots,m-1,
$$
\noindent
with $m$ such that $N=2^m$.

This method will be also considered for the Manneville-Pomeau
processes when $s\geq 1$. This corresponds to the case when the
invariant measure $\mu_s$ is not a probability measure (see Table
$7.1$).

\vspace{0.3cm}

\section{\normalsize Monte Carlo Simulation for the ``Long Dependence" Case}
\renewcommand{\theequation}{5.\arabic{equation}}
\setcounter{equation}{0}

In this section we present the Monte Carlo simulation results
comparing the six different estimation methods given in Section $4$
for the {\it long dependence case}.

Let $\{X_t\}_{t\in \N}$ be the Manneville-Pomeau process, given by
the expression (\ref{MPprocess}), where $\varphi=\I_A$ with
$A={(0.1,0.9)}$ such that $X_t= \I_{A}\circ T_{s}^t$.

One chooses at random a value $x_0$ of the random
variable $X_0$ according to a uniform distribution
(this is the same as to choose $x_0$ at random according
to the probability $\mu_s$).
Let $\{X_t\}_{t=0}^{N-1}$ be a time series with $N$
observations from the process $\{X_t\}_{t\in \N}$
obtained from such $x_0$.
Hence, this time series is given by
\begin{equation}\label{MPprocess2}
X_t= \I_{A}(T_s^t (x_0))=
\I_{(0.1,0.9)}(T^t_s(x_0)), {\mbox { for all }}
t=0,\cdots ,N-1.
\end{equation}

\newpage


\vspace{0.4cm}

The simulations presented here are based on such time series.

Figures $5.1$ (a) and (b) present the sample autocorrelation and the
periodogram functions, respectively, for a time series with sample
size $N=10,000$ obtained from (\ref{MPprocess2}) when $s=0.8$.

The following results were obtained from Monte Carlo simulations in
Fortran routines and using the IMSL library. We remark that for the
long dependence case one needs large number of sets of data
requiring high computational time.

\vspace{0.3cm}


\noindent
\begin{center}
{\bf Table 5.1}:  Estimation results when
$s\in \{0.60,0.65\}$.
\end{center}

\begin{center}\label{table51}
{\scriptsize
\begin{tabular}{|c|c|c|c|c|c|c|}

\hline\hline
$s$&$N$&{\it Method}&{\it mean}$({\hat s})$
&{$sd(\hat s)$}&${mse(\hat s)}$\\
\hline\hline

&&$Perio$&0.6545&0.1394&0.0223\\
&&$Parzen$&0.6313&0.1125& 0.0136\\
&10,000&$Cos$(1)&0.5531&0.0572&0.0054\\
&&$Cos$(2)&0.5993&0.0220&{\bf 0.0005}\\
&&$Varmp$&0.5309&0.0396& 0.0063\\
&&$Vpmp$&0.5598&0.0718&0.0067\\
\cline{2-6}

&&$Perio$&0.6364&0.1094&0.0130\\
&&$Parzen$&0.6147&0.0086&0.0070\\
&20,000&$Cos$(1)&0.5488&0.0535&0.0054\\
0.60&&$Cos$(2)&0.5979&0.0264&{\bf 0.0007}\\
&&$Varmp$&0.5241&0.0303&0.0067\\
&&$Vpmp$&0.5513&0.0583&0.0057\\
\cline{2-6}

&&$Perio$&0.6004&0.1051&0.0110\\
&&$Parzen$&0.5865&0.0736&0.0056\\
&30,000&$Cos$(1)&0.5275&0.0508&0.0078\\
&&$Cos$(2)&0.5933&0.0316&{\bf 0.0010}\\
&&$Varmp$&0.5204&0.0264&0.0070\\
&&$Vpmp$&0.5144&0.0608&0.0110\\
\hline\hline

&&$Perio$&0.7539&0.1518&0.0337\\
&&$Parzen$&0.7107&0.0107&0.0151\\
&10,000&$Cos$(1)&0.6145&0.0614&0.0050\\
&&$Cos$(2)&0,6129&0,0198&{\bf 0,0017}\\
&&$Varmp$&0.5293&0.0332&0.0156\\
&&$Vpmp$&0.5461&0.0763&0.0166\\
\cline{2-6}

&&$Perio$&0.7113&0.0779&0.0098\\
&&$Parzen$&0.6927&0.0706&0.0068\\
&20,000&$Cos$(1)&0.6035&0.0472&0.0044\\
0.65&&$Cos$(2)&0.6076&0.0181&{\bf 0.0021}\\
&&$Varmp$&0.5251&0.0246&0.0162\\
&&$Vpmp$&0.5257&0.0630&0.0194\\
\cline{2-6}

&&$Perio$&0.6806&0.0445&0.0029\\
&&$Parzen$&0.6910&0.0392&0.0032\\
&30,000&$Cos$(1)&0.6141&0.0552&0.0043\\
&&$Cos$(2)&0.6090&0.0147&{\bf 0.0019}\\
&&$Varmp$&0.5419&0.0262&0.0123\\
&&$Vpmp$&0.5451&0.0562&0.0141\\
\hline\hline
\end{tabular}
}
\end{center}

\vspace{0.3cm}

For all tables presented here, we calculated the mean ($mean$), the
standard deviation ($sd$) and the mean squared error ($mse$) values
for all estimators of $s$. The smallest mean squared error is shown
in boldfaced character in these tables. All simulations are based in
$200$ replications unless for Tables $5.3$ and $5.4$ where we use
$50$ replications. For the estimator $Cos$ we used two different
values for the limit $g(N)=N^{\alpha_i}$: $Cos(1)$ means $\alpha_1
=0.5$ and $Cos(2)$ means $\alpha_2 =0.7$.

Table $5.1$ presents the results for the six estimation methods
proposed in Section $4$ for the {\it long dependence case\/} for
$s\in \{0.60,0.65\}$ and for three different values of $N$.

From Table $5.1$ we observe that the estimators $Varmp$ and $Vpmp$
are very much biased: this was also true for the ARFIMA process (see
Olbermann, 2002). The best result, in terms of small mean squared
error value is the estimator $Cos(2)$ for both values of $s$ and for
any sample size considered.

In Table $5.2$ we present the results for the case when $s=0.80$
considering only the sample size $N=10,000$. The best results were
for the methods $Perio$ and $Parzen$ since the other methods have
higher bias. As $s$ approaches to the value $1$, the time series
$\{X_t\}_{t=0}^{N-1}$, given by (\ref{MPprocess2}), stays long time
in zero, resulting in very poor estimates. The methods $Varmp$ and
$Vpmp$ are not recommended in this situation due to their higher
bias.

\vspace{0.4cm}

\noindent
\begin{center}
{\bf Table 5.2}:  Estimation Results when $s=0,80$.
\end{center}

\begin{center}\label{table52}

\begin{tabular}{|c|c|c|c|c|c|}
\hline\hline
$N$&{\it Method}&{\it mean}$({\hat s})$
&{$sd(\hat s)$}&${mse(\hat s)}$\\
\hline\hline

&$Perio$&0.7773&0.1648&0.0275\\
&$Parzen$&0.7607&0.1444&{\bf 0.0222}\\
10,000&$Cos$(1)&0.6286&0.2507&0.0919\\
&$Cos$(2)&0.6626&0.0822&0.0256\\
&$Varmp$&0.5472&0.0426&0.0657\\
&Vpmp&0.5781&0.0806&0.0557\\
\hline\hline

&$Perio$&0.6921&0.1220&{\bf 0.0264}\\
&$Parzen$&0.6740&0.1127&0.2849\\
20,000&$Cos$(1)&0.5731&0.0699&0.0563\\
&$Cos$(2)&0.6434&0.0437&{\bf 0.0264}\\
&$Varmp$&0.5292&0.0434&0.0752\\
&$Vpmp$&0.5416&0.0848&0.0739\\
\hline\hline

&$Perio$&0.6559&0.1164&0.0342\\
&$Parzen$&0.6150&0.1044&0.0456\\
30,000&$Cos$(1)&0.5335&0.1713&0.1002\\
&$Cos$(2)&0.6382&0.0337&{\bf 0.0273}\\
&$Varmp$&0.5354&0.0524&0.0727\\
&$Vpmp$&0.5434&0.0861&0.0732\\
\hline\hline
\end{tabular}
\end{center}

\vspace{0.4cm}

The simulations presented in Tables $5.3$ and $5.4$ are based on
$50$ replications.

\vspace{0.2cm}

In Table $5.3$ we present the results for the wavelet method. We
only consider the bases Haar and Mexican hat. We remark that this
estimator requires power of two for the sample size. This table
presents the results when $s\in \{0.65,0.80\}$ with three different
values for $N$. We observe that the Mexican hat basis has advantages
over the Haar basis presenting smaller bias and mean squared error
values.

\vspace{0.2cm}

After the analysis of the {\it long dependence case\/} we make a few
comments about another regime, that is, when $s\geq 1$.


\vspace{5.5cm}

\noindent
\begin{center}
{\bf Table 5.3}:  Estimation results when $s\in \{0.65,0.80\}$.
\end{center}

\begin{center}\label{table53}

\begin{tabular}{|c|c|c|c|c|c|c|}
\hline\hline $s$&$N$&{\it Wavelet Basis}&{\it mean}$({\hat s})$
&{$sd(\hat s)$}&${mse(\hat s)}$\\
\hline\hline

&8,192&Haar&0.8531&0.0470&0.0434\\
&&Mexican hat&0.8022&0.0480&{\bf 0.0254}\\
\cline{2-6}

&16,384&Haar&0.8311&0.0446&0.0347\\
0.65&&Mexican hat&0.7882&0.0472&{\bf 0.0213}\\
\cline{2-6}

&32,768&Haar&0.8283&0.0619&0.0355\\
&&Mexican hat&0.7864&0.0451&{\bf 0.0206}\\
\hline\hline

&8,192&Haar&0.9839&0.0619&0.0376\\
&&Mexican hat&0.8873&0.0670&{\bf 0.0120}\\
\cline{2-6}

&16,384&Haar&0.9321&0.0659&0.0217\\
0.80&&Mexican hat&0.8237&0.0675&{\bf 0.0050}\\
\cline{2-6}

&32,768&Haar&0.8639&0.0915&0.0120\\
&&Mexican hat&0.7747&0.0464&{\bf 0.0027}\\
\hline\hline
\end{tabular}
\end{center}

\vspace{0.3cm}

In Table $5.4$ we present the case where $s\geq 1$ meaning that the
invariant measure $\mu_s$ does not correspond to a probability
measure for the process $\{X_t\}_{t\in \N}$, given by
(\ref{MPprocess}). This table presents values of
$s\in\{1.0,1.1,1.2,1.3\}$ and sample size $N=32,768$. The best
results were for the Haar basis. Notice that when $s\geq 1$ any
method based on the periodogram function does not make sense (for
the process obtained from the iterations of the Manneville-Pomeau
transformation $T_s$ when $x_0$ is chosen at random).

An interesting question to be investigated: is it true that for any
deterministic (such as Manneville-Pomeau, Infinite Markov Chain,
etc$\ldots$) or purely stochastic process (such as ARFIMA,
etc$\ldots$) depending only on the decay of the rate of convergence
of the autocorrelation function, there exists a better wavelet basis
(such as Haar, Mexican hat, Shannon, etc$\ldots$) to estimate the
exponent of decay?

\vspace{0.3cm}

\noindent
\begin{center}
{\bf Table 5.4}:  Estimation results when $s\geq 1$ with $N=32,768$.
\end{center}

\begin{center}\label{table54}

\begin{tabular}{|c|c|c|c|c|c|}
\hline\hline $s$&{\it Wavelet Basis}&{\it mean}$({\hat s})$
&{$sd(\hat s)$}&${mse(\hat s)}$\\
\hline\hline

1.0&Haar&0.9461&0,1090&{\bf 0.0145}\\
&Mexican hat&0.8931&0.1148&0.0243\\
\hline

1.1&Haar&1.0924&0.0589&{\bf 0.0034}\\
&Mexican hat&0.9943&0.0461&0.0132\\
\hline

1.2&Haar&1.0825&0.0729&{\bf 0.0190}\\
&Mexican hat&0.9642&0.0939&0.0642\\
\hline

1.3&Haar&1.1422&0.0638&{\bf 0.0288}\\
&Mexican hat&1.0064&0.0703&0.0910\\
\hline\hline
\end{tabular}
\end{center}

\vspace{0.3cm}

\section{\normalsize Estimation in the ``Not So Long Dependence" Case}
\renewcommand{\theequation}{6.\arabic{equation}}
\setcounter{equation}{0} \vspace{0.5cm}

In the {\it not so long dependence case\/} one can estimate the
value $s$ using the exactly $a$-H\"older property in the point
$x_0=0$ (see Bary, 1964 and Fisher and Lopes, 2001). Suppose
$$
a\approx \frac {\ln (|f_X(x_0) - f_X(y)|)}{\ln (|x_0 - y|)}, {\mbox
{ for }} y \in (-\pi, \pi) {\mbox { very close to zero, }}
$$
where $f_X(\cdot)$ is the spectral density function, given in
(\ref{spectral}), of the process $\{X_t\}_{t\in \N}$ given in
(\ref{MPprocess}). We then define the estimator
\begin{equation}\label{holder}
\hat s=\frac {1}{a+2} \,\, {\mbox { where }} \,\, a=\frac {\ln
(|I(\omega_0) - I(\omega_j)|)} {\ln (|\omega_0 - \omega_j|)},
\end{equation}
\noindent with $I(\cdot)$ the periodogram function, given by
(\ref{periodogram}), with $\omega_0=0$ and $\omega_j$ is a Fourier
frequency, given by (\ref{omegah}), very close to zero.

The main goal of this section is to describe two different
estimation methods to estimate the transformation $T_s$, or
equivalently, to estimate the parameter $s$, when $s\in
(0,\frac12)$. For this purpose, we consider a finite time series
$\{X_t\}_{t=0}^{N-1}$ obtained from the process $\{X_t\}_{t\in\N}$
given by (\ref{MPprocess2}). The two methods are proposed by
(\ref{holder}) when the periodogram or its smoothed version by the
Parzen lag window functions are used.

These methods are described in this section and in Section $7$ we
present a Monte Carlo simulation study comparing them.

\vspace{1cm}

\noindent {\bf {\it P} Estimator }

\vspace{.5cm}

This estimation method is based on the expression (\ref{holder})
above where $I(\cdot)$ is the periodogram function given by the
expression (\ref{periodogram}). We denote it by {\it P}.

\vspace{0.5cm}

\noindent {\bf {\it SP} Estimator }

\vspace{.5cm}

This estimation method is based on the expression (\ref{holder})
above where the periodogram function $I(\cdot)$ is now replaced by
the smoothed periodogram function $f_{sm}(\cdot)$ using the Parzen
spectral window. We denote this estimator by {\it SP}.

\vspace{0.3cm}

\section{\normalsize Monte Carlo Simulation for the ``Not So Long Dependence" Case}
\renewcommand{\theequation}{7.\arabic{equation}}
\setcounter{equation}{0}

In this section we present the Monte Carlo simulation results
comparing the two methods given in Section $6$ for the {\it not so
long dependence case}.

Let $\{X_t\}_{t\in \N}$ be the Manneville-Pomeau process, given by
the expression (\ref{MPprocess2}), where $\varphi=\I_A$ with
$A={(0.1,0.9)}$ such that $X_t= \I_{A}\circ T_{s}^t$.

One chooses at random a value $x_0$ of the random variable $X_0$
according to a uniform distribution (this is the same as to choose
$x_0$ at random according to the probability $\mu_s$). Let
$\{X_t\}_{t=0}^{N-1}$ be a time series with $N$ observations
obtained from (\ref{MPprocess2}). The simulations presented here are
based on such time series and were obtained by Fortran routines with
some help of the IMSL library.

In Table $7.1$ we present some simulation results for the {\it not
so long dependence case} based on the two methods reported in
Section $6$. We calculated the mean ($mean$), the standard deviation
($sd$) and the mean squared error ($mse$) values for each method.
The smallest mean squared error is shown in boldfaced character in
this table. These simulations are based in 200 replications with
$s\in \{0.35,0.40,0.45\}$ and two different sample sizes $N$. Note
that as we have a better mixing rate of convergence for the {\it not
so long dependence case\/} the biases here are smaller than the case
of long dependence.

\vspace{3.3cm}

\noindent
\begin{center}
{\bf Table 7.1}:  Estimation results for $s\in \{0.35,0.40,0.45\}$.
\end{center}

\begin{center}\label{table7}

\begin{tabular}{|c|c|c|c|c|}
\hline\hline
$s$&$N$&Statistics&$P$&$SP$\\
\hline\hline
&&{\it mean}$({\hat s})$&0.4078&0.3970\\
&10,000&{$sd(\hat s)$}&0.0374&0.0255\\
0.35&&${mse(\hat s)}$&0.0047&{\bf 0.0028}\\
\cline{2-5}
&&{\it mean}$({\hat s})$&0.3870&0.4136\\
&30,000&{$sd(\hat s)$}&0.0298&0.0208\\
&&${mse(\hat s)}$&{\bf 0.0022}&0.0044\\
\hline\hline

&&{\it mean}$({\hat s})$&0.4210&0.4024\\
&10,000&{$sd(\hat s)$}&0.0378&0.0258\\
0.40&&${mse(\hat s)}$&0.0018&{\bf 0.0006}\\
\cline{2-5}
&&{\it mean}$({\hat s})$&0.4397&0.4046\\
&30,000&{$sd(\hat s)$}&0.0432&0.0405\\
&&${mse(\hat s)}$&0.0034&{\bf 0.0016}\\
\hline\hline

&&{\it mean}$({\hat s})$&0.4652&0.4359\\
&10,000&{$sd(\hat s)$}&0.0312&0.0285\\
0.45&&${mse(\hat s)}$&{\bf 0.0012}&0.0050\\
\cline{2-5}
&&{\it mean}$({\hat s})$&0.5218&0.4808\\
&30,000&{$sd(\hat s)$}&0.0800&0.0619\\
&&${mse(\hat s)}$&0.0115&{\bf 0.0047}\\
\hline\hline
\end{tabular}
\end{center}

\vspace{0.3cm}

\section{\normalsize Conclusions}
\renewcommand{\theequation}{8.\arabic{equation}}
\setcounter{equation}{0}

We analyzed the estimation of the parameter $s$ in the
Manneville-Pomeau processes in {\it the long} and {\it not so
long-range dependence} cases.

We considered several estimation methods for both situations and we
compared them with the method called here {\it Perio} presented by
Schuster (1984) and largely used by the physicists. We point out
that from the analysis of the simulation studies in Tables $5.1$ to
$5.3$ the wavelet method ($Wmp$) had a better performance than the
{\it Perio\/} method. One can see this in the case when $s=0.8$
where the estimator $Wmp$ from the Mexican hat wavelet basis gave
the best result in terms of small mean squared error.

The methods $Varmp$ and $Vpmp$ presented
the higher biases while the method $Cos(2)$
had the best results for the cases when
$s\in \{0.60,0.65\}$, with the smallest
mean squared error for the sample size
analyzed.

We analyzed the performance of the
method $Wmp$ based on the wavelet theory
for the Manneville-Pomeau processes
when $s\geq 1$, that corresponds to the
situation where the invariant measure
$\mu_s$ is not a probability measure.
In this case, the best results were
obtained when the Haar basis was considered.

Among the estimation methods proposed for the {\it not so long
dependence case\/} the one based on the smoothed periodogram
function using the Parzen spectral window had the best results with
lower bias and mean squared error values. The wavelet method $Wmp$
also works fine in this case, but we do not present here simulation
results for this case.

\vspace{0.3cm}

\subsection*{Appendix A}
\renewcommand{\theequation}{A.\arabic{equation}}
\setcounter{equation}{0}

\vspace{0.5cm}

Let $\{X_t\}_{t\in \N}$ be the Manneville-Pomeau process defined in
(\ref{MPprocess}). Let $\rho_X(\cdot)$ and ${f_X(\cdot)}$ be,
respectively, the autocorrelation and the spectral density functions
of this process.

In this appendix we present some general properties
of the Fourier series. In this way we will explain
why the hyperbolic (or polynomial) decay of the
autocorrelation function, that is,

$$
{\rho_X(h)}\approx {h^{-u}}, \, \mbox{ for  }
0<\mu<1,
$$
corresponds to

$$
{f_X(\lambda)}\approx {\lambda^{u-1}},
$$
for the spectral density function of the process $\{X_t\}_{t\in \N}$
given by (\ref{MPprocess}).

First we will explain the {\it not so long
dependence case\/}.

If the function $g$ is $n$-times differentiable
and $g^n(\cdot)$ is $a$-H\"older with
$0<a<1$, we say that $g$ is $(n+a)$-H\"older.

The relationship of the hyperbolic decay between the autocorrelation
function of the Manneville-Pomeau process and its spectral density
function is only a question related to Fourier series (see Bary,
1964).

\vspace{0.5cm}

\noindent
{\bf Theorem A.1:} {\it Suppose that
$b_n\approx n^{-u}$, for some $u$, and suppose also that
$g(\theta)= \sum_{n=1}^\infty   b_n \cos (n \, \theta)$,
where $b_n\in {\R}$, converges to zero. If $a$ is positive
and $g(\cdot)$ is a H\"older function of order $a$,
then there exists a positive constant $c$ such that
$b_n < c \, n^{-(1+a)}$, for all $n\in \N -\{0\}$.}

\vspace{0.5cm}

\noindent
{\bf Theorem A.2:} {\it Suppose that
$b_n\approx n^{-u}$, for some $u$, and suppose also that
$g(\theta)= \sum_{n=1}^\infty   b_n \cos (n \, \theta)$,
where $b_n \in {\R}$, decreases monotonously to zero.
If $a$ is positive and there exists a positive constant
$c$ such that $b_n < c\, n^{-(1+a)}$,
then $g(\cdot)$ is a H\"older function of order $a$.}

\vspace{0.3cm}

Theorems A.1 e A.2 (see Chapter II, Section $3$ and Chapter X,
Section $9$, respectively, in Bary, 1964) apply to the {\it not so
long dependence case\/}.

Another interesting result on Fourier series,
that can be applied now for the
{\it long dependence case\/}, is described in
the next theorem.

\vspace{0.5cm}

\noindent
{\bf Theorem A.3 (Riesz):} {\it Suppose that
$g(\theta)= \sum_{n=1}^\infty   b_n
\cos (n \, \theta)$, for all $\theta\in (-\pi,\pi)$
and that $b_n \in {\R}$ is such that the sequence
$\{b_n\}_{n\in\N}$ decreases monotonously
to zero when $n\to \infty$.
Suppose there exists a positive real constant
$u$ such that
$b_n\approx n^{-u}$.
Suppose there exists also a positive real constant
$b \in (-1,0)$ such that
$$
|g(\theta)|\approx  |\theta|^{b}.
$$

\vspace{0.2cm}

\begin{itemize}
\item [$(a)$] If there exist $a\in (-1,0)$, $\epsilon>0$ and a
positive real constant $k$ such that
$$
\left |\frac{g(\theta)}{\theta^{a} }\right| \leq k,
\, \mbox{ for all } \, 0<\theta<\epsilon,
$$
then $u\geq 1+a$. That is, the decreasing
velocity of
$|b_n|$ is at least of order $ n^{-(1+a)}$,
when $n\to \infty$.

\vspace{0.2cm}

\item[$(b)$] If there exist $a\in (-1,0)$ and a positive real constant
$v$ such that $|b_n|< v\, n^{-(1+a)}$, then $b\leq a$. That is,
$g(\theta)$ is at least of order of $|\theta|^a$, when $\theta\to
0$.
\end{itemize}

Hence, from $(a)$ and $(b)$ above one concludes that $u=1+b$.}

\vspace{0.4cm}

\noindent {\bf Remark A.1:} In the general cases, we point out that
there exist sequences $\{b_n\}_{n\in\N}$ (not monotonous) such that
$c_1 n^{-u}<  |b_n| < c_2 n^{-u}$, for some positive constants $c_1$
and $c_2$ and $u$ such that $0<u<1$, but $g(\theta)$ does not
satisfy $ c_3 |\theta|^{b}\leq |g(\theta)|\leq c_4 |\theta|^{b}$ for
any fixed positive constants $c_3, c_4$ and $b$.

\vspace{0.4cm}

Theorem A.$3$ is a consequence of the following result.

\vspace{0.5cm}

\noindent {\bf Theorem A.4 (Riesz):} {\it Suppose that $g(\theta)=
\sum_{n=1}^\infty b_n \cos (n \, \theta)$, for all $\theta\in
(-\pi,\pi)$, and that $\{b_n\}_{n\in\N} \in {\R}$ decreases
monotonously to zero. Let $p>1$ and $q>1$ be such that $ \frac{1}{p}
+ \frac{1}{q}=1$.}

\vspace{0.2cm}
\begin{itemize}
\item[$(a)$] {\it If $g \in {\cal L}^p$, then $ \sum_{n=1}^\infty   |b_n|^q <
\infty$.}

\vspace{0.2cm}

\item[$(b)$] {\it If $ \sum_{n=1}^\infty |b_n|^q < \infty$, then $g \in {\cal
L}^p$.}
\end{itemize}

\vspace{0.4cm}

\noindent {\bf Remark A.2:} Theorem A.$3$ follows from Theorem A.$4$
making use of

\vspace{0.2cm}

\noindent
\begin{itemize}
\item[(a)] for any continuous function $f$, on $(0,\pi)$, of order
$x^\alpha$ ($x$ close to zero), then $f\in {\cal L}^1\Leftrightarrow
\alpha> -1$ \noindent and

\item[(b)] for any sequence $c_n$ of order $n^{-\beta}$ ($n$ close
to infinity) then $ \sum_{n=1}^\infty|c_n| < \infty \Leftrightarrow
\beta> 1$.
\end{itemize}

\vspace{0.4cm}

Theorem A.$4$ follows easily from the first theorem of Chapter X,
Section $9$ of Bary (1964).

The above results justify the ideas used in the estimation methods
$Perio$, $Parzen$, $Cos(1)$ and $Cos(2)$, given in Section $4$.

\vspace{0.3cm}

\subsection*{Appendix B}
\renewcommand{\theequation}{B.\arabic{equation}}
\setcounter{equation}{0}

\vspace{0.5cm}

Considering the rate of convergence to zero of the autocorrelation
function one can also get an estimate of the order of magnitude of
the variance for the partial sums $S_N=\sum^{N-1}_{i=0} X_i$ from a
time series $\cdots X_{-3}, X_{-2}, X_{-1},X_0, X_1, \cdots ,
X_{N-1}$. In Proposition B.$1$ below we present a proof of the
estimated value for the variance of the random variable $S_N$. In
Proposition B.$2$ we give a precise estimate of the order of growth
for the variance of this random process.

We point out that the stationary process stated above and given by
$$
X_t=(\varphi\circ T_s^t)(X_0), \mbox{ for } t\in {\N},
$$
can be considered defined for all $t\in \Z$, via the natural
extension transformation (see section $5.3$ in Lopes and Lopes,
1998).

\vspace{.5cm}

\noindent
{\bf Proposition B.1}: {\it Let $\{X_t\}_{t\in{ \Z}}$
be any stationary stochastic process.
Let $S_N=\sum^{N-1}_{i=0} X_i$ be the partial sum
of a time series $X_0, X_1, \cdots , X_{N-1}$ from this
process. Then,

$$
Var(S_N)= 2N \, \left[\frac{\gamma_X(0)}{2} +
\frac{1}{N} \sum^{N-1}_{j=1} (N-j)\,
  \gamma_X(j)\right],
$$
\noindent
where $\gamma_X(\cdot)$ is the autocovariance function
of the process $\{X_t\}_{t\in{ \Z}}$.}

\vspace{0.5cm}

\noindent
{\bf Proof:} Since the process $\{X_t\}_{t\in{ \Z}}$
is stationary, we observe that
\begin{eqnarray}\label{B1}
Var(S_N) &=& Var\left(\sum^{N-1}_{i=0} X_i\right) \nonumber \\
&=& \sum^{N-1}_{i=0} Var (X_i) + \sum^{N-1}_{j=0}
\sum^{N-1}_{\ell=0} cov (X_j, X_\ell)\nonumber \\
&=& N Var(X_0)+ \sum^{N-1}_{j=0} \sum^{N-1}_{\ell=0}
\left(\E (X_j X_\ell) -
 [\E(X_0)]^2 \right)\nonumber \\
& = & N \, \gamma_X(0)+  2 \sum^{N-1}_{\stackrel{j,l=0}
{j<\ell}} \gamma_X(j-\ell ).
\end{eqnarray}

\noindent It follows from the expression (\ref{B1}) that

\begin{eqnarray}\label{B2}
\nonumber
Var(S_N) &=& N \, \gamma_X(0)+
 2 \sum^{N-1}_{\stackrel{j,l=0}{j<\ell}} \gamma_X(j-\ell)\\
 \nonumber
 &=& N \, \gamma_X(0)+2\left( \underbrace{\gamma_X(-1) +\gamma_X(-2) +
 \gamma_X(-3)+
 \cdots + \gamma_X(-N+1)}_{j=0} \right.\\
\nonumber
&+&  \underbrace{
\gamma_X(-1) +\gamma_X(-2) + \cdots + \gamma_X(1-(N-1))}_{j=1}\\
\nonumber &+&  \underbrace{ \gamma_X(-1) +\gamma_X(-2) + \cdots +
\gamma_X(2-(N-1))}_{j=2}\\
\nonumber &+& \left. \underbrace{ \gamma_X(-1) +\gamma_X(-2) +
\cdots + \gamma_X(3-(N-1))}_{j=3} + \cdots + \underbrace{
\gamma_X(-1) }_{j=N-2} \right)\\
\nonumber & = & N\,  \gamma_X(0) + 2 \left[ (N-1) \gamma_X(-1) +
(N-2) \gamma_X(-2) + (N-3) \gamma_X(-3)  \right.\\
\nonumber &+& \left. \cdots + 3\gamma_X(-(N-3)) + 2\gamma_X(-(N-2))
+ \gamma_X(-(N-1))\right]\\
&=& N \, \gamma_X(0) + 2 \sum^{N-1}_{j=1} (N-j)\, \gamma_X(-j)=
 N \, \gamma_X(0) + 2 \sum^{N-1}_{j=1} (N-j)\,
\gamma_X(j).
\end{eqnarray}

\noindent
The last equality (\ref{B2}) follows from the fact that
the process is stationary. This implies that
$\gamma_X(j)=\gamma_X(-j)$.

Therefore,
$$
Var(S_N)= N \, \gamma_X(0) + 2 \sum^{N-1}_{j=1}
(N-j)\,  \gamma_X(j),
$$
\noindent and this completes the proof of Proposition B.$1$. \, \qed

\vspace{0.3cm}

In the next proposition we show the order
of $Var(S_N)$, with respect to $N$,
for a quite general class of stationary stochastic processes.

\vspace{0.5cm}

 \noindent
{\bf Proposition B.2}: {\it Let $\{X_t\}_{t\in{ \Z}}$ be
any stationary stochastic process. Let
$S_N=\sum_{i=0}^{N-1}X_i$ be the partial sum of a time
series $X_0,X_1,\cdots,X_{N-1}$ from the process
$\{X_t\}_{t\in{ \Z}}$. If there exists $u\in (0,1)$
such that $\gamma_X(h)\approx h^{-u}$, then
$$
Var(S_N)\approx N^{2 - u}.
$$
}

\vspace{0.2cm}

\noindent
{\bf Proof:} For $u\in (0,1)$, the integral
$$
I=\int_0^1 (1-x) x^{-u} dx
$$
\noindent
is finite. Then, for any $N\in {\N}$,
one can consider the Riemann sums
associated to the partition
$$
\left\{0, \frac{1}{N}, \frac{2}{N},\cdots, \frac{N-1}{N},
1\right\},
$$
\noindent
obtaining the approximation

$$
\sum_{j=1}^N \left(1-\frac{j}{N}\right)
\left(\frac{j}{N}\right)^{-u} \frac{1}{N},
$$
\noindent
that converges to $I$, when $N\rightarrow \infty$.
\newline

From similar arguments proposed in Lemma $8.1$ of Fisher and Lopes
(2001),  consider
\begin{eqnarray}\label{B3}
c_N &=& \sum_{j=1}^N \left(1-\frac{j}{N}\right)
\left(\frac{j}{N}\right)^{-u} \nonumber \\
&=& \sum_{j=1}^N \left(\frac{N-j}{N}\right)
\left(\frac{j}{N}\right)^{-u} \nonumber \\
&=&\sum_{j=1}^N (N-j) \, j^{-u}
\left(\frac{1}{N}\right)^{1-u}.
\end{eqnarray}
\noindent
Given $\varepsilon >0$, for $N$ sufficiently
large, one has that

$$
I-\varepsilon\,  \leq \frac{1}{N} \, c_N \,\leq
 I+ \varepsilon.
$$
\noindent Using the expression (\ref{B3}), the above inequality is
given by

\begin{equation}\label{B4}
 (I-\varepsilon)\, N^{1-u} \leq\frac{1}{N}
\sum_{j=1}^N (N-j) \, j^{-u}\leq
(I+ \varepsilon) \, N^{1-u},
\end{equation}
\noindent
for $N$ sufficiently large.

Therefore,

$$
\frac{1}{N}\sum_{j=1}^N (N-j) \, j^{-u} \,
\mbox{ is of order } N^{1-u}.
$$
\newline

\noindent From the expressions (\ref{B2}) and (\ref{B4}) one has
\begin{equation*}
  Var(S_N)= 2N \, \left[\frac{\gamma_X(0)}{2} +
\frac{1}{N} \sum^{N-1}_{j=1} (N-j)\
  \gamma_X(j)\right]\approx N^{-u},
\end{equation*}
\noindent and this completes the proof of Proposition B.$2$. \, \qed

\vspace{0.2cm}

The above results justify the ideas used in the estimation methods
$Varmp$ and $Vpmp$, given in Section $4$.

\vspace{0.3cm}
\subsection*{Acknowledgements}

\vspace{0.3cm}

\noindent S.R.C. Lopes was partially supported by CNPq-Brazil, by
Pronex {\it Probabilidade e Processos Estoc\'asticos} (Conv\^enio
MCT/CNPq/FAPERJ - Edital 2003), by Edital Universal {\it Modelos com
Depend\^encia de Longo Alcance: An\'alise Probabil\'istica e
Infer\^encia} (CNPq-No. 476781/2004-3) and also by {\it Funda\c
c\~ao de Amparo \`a Pesquisa no Estado do Rio Grande do Sul\/}
(FAPERGS Foundation). A.O. Lopes was partially supported by
CNPq-Brazil,  by Pronex {\it Sistemas Din\^amicos}, by the {\it
Millennium Institute in Mathematics\/} and is also beneficiary of
CAPES financial support.

\vspace{0.5cm}


\noindent {\bf REFERENCES \rm}

\vspace{0.5cm}

{\footnotesize

\begin{description}

\item Absil, P.A., R. Sepulchre, A. Bilge and P. G\'erard (1999).
``Nonlinear Analysis of Cardiac Rhythm Fluctuations using DFA
Method". {\it Physica A\/}, {\bf 272}, pp. 235-244. \vspace*{0.1cm}

\item Bary, N.K. (1964). {\it A  Treatise on Trigonometric series},
Vol. I and II. New York: Pergamon Press. \vspace*{0.1cm}

\item Beran, J. (1994). {\it Statistics for Long-Memory Processes}.
New York: Chapman \& Hall. \vspace*{0.1cm}

\item Brockwell, P.J. and  R.A. Davis (1991). {\it Time Series:
Theory and Methods}. New York: Springer-Verlag. \vspace*{0.1cm}

\item Chazottes, J.-R., E. Floriani and R. Lima (1998). ``Relative
Entropy and Identification of Gibbs Measures in Dynamical Systems".
{\it Journal of Statistical Physics}, {\bf 90}(3-4), pp. 697-725.
\vspace*{0.1cm}

\item Chazottes, J.-R., P. Collet and B. Schmitt (2005).
``Statistical Consequences of the Devroye Inequality for Processes.
Applications to a Class of Non-Uniformly Hyperbolic Dynamical
Systems". {\it Nonlinearity}, {\bf 18}(5), pp. 2341-2364.
\vspace*{0.1cm}

\item Collet, P., A. Galves and A. Lopes (1995). ``Maximum
Likelihood and Minimum Entropy Identification of Grammars". {\it
Random Computation and Dynamics}, {\bf 3}(4), pp. 241-250.

\item Collet, P., S. Martinez and B. Schmitt (2004). ``Asymptotic Distribution
of Tests for Expanding Maps in the Interval". {\it Ergodic Theory
and Dynamical Systems}, {\bf 24}(3), pp. 707-722. \vspace*{0.1cm}

\item Collet, P. (2005). {\it Dynamical Systems and Stochastic
Processes}. Notas de Los Cursos, XIV Escuela Latinoamericana de
Matem\'atica. Montevideo: Facultad de Ciencias-Facultad de
Ingenieria, pp. 51-186. \vspace{0.1cm}

\item Collet, P. and J.-P. Eckmann (2006). {\it Concepts and Results
in Chaotic Dynamics: a Short Course}. Theoretical and Mathematical
Physics. Berlin: Springer-Verlag. \vspace*{0.1cm}

\item Feller, W. (1949). ``Fluctuation Theory of Recurrent
Events''. {\it Transactions of the American Mathematical Society},
{\bf 67}, pp. 98-119. \vspace*{0.1cm}

\item Fisher, A. and A. Lopes (2001). ``Exact bounds for the
polynomial decay of correlation, 1/f noise and the CLT for the
equilibrium state of a non-H\"older potential''.  {\it
Nonlinearity}, {\bf 14}, pp. 1071-1104. \vspace*{0.1cm}

\item Geweke, J. and S. Porter-Hudak (1983). ``The Estimation and
Application of Long Memory Time Series Model''. {\it Journal of Time
Series Analysis}, {\bf 4}(4), pp. 221-238. \vspace*{0.1cm}

\item Gou\"ezel, S. (2004). ``Central Limit Theorem and Stable Laws
for Intermittent Maps". {\it Probability Theory and Related Fields},
{\bf 128}(1), pp. 82-122. \vspace*{0.1cm}

\item Guharay, S., B.R. Hunt, J.A. Yorke and O.R. White (2000).
``Correlations in DNA Sequences Across the Three Domains of Life".
{\it Physica D}, {\bf 146}(1-4), pp. 388-396. \vspace*{0.1cm}

\item Jensen, M.J. (1999). ``Using Wavelets to Obtain a Consistent
Ordinary Least Square Estimator''. {\it Journal of Forecasting},
{\bf 18}, pp. 17-32. \vspace*{0.1cm}

\item Lopes, A. (1993).
 ``The zeta function, non-differentiability of
pressure and the critical exponent of transition''. {\it Advances in
Mathematics}, {\bf 101}, pp. 133-165. \vspace*{0.1cm}

\item Lopes, A. and S.R.C. Lopes (1998). ``Parametric Estimation
and Spectral Analysis of Piecewise Linear Maps of the Interval".
{\it Advances in Applied Probability\/}, {\bf 30}(3), pp. 757-776.
\vspace*{0.1cm}

\item Lopes, A. and S.R.C. Lopes (2002). ``Convergence in
Distribution of the Periodogram for Chaotic Processes". {\it
Stochastics and Dynamics\/}, {\bf 2}(4), pp. 609-624.
\vspace*{0.1cm}

\item Lopes, S.R.C., B.P. Olbermann and V.A. Reisen (2002).
``Non-stationary Gaussian ARFIMA Processes: Estimation and
Application". {\it Brazilian Review of Econometric\/}, {\bf 22}(1),
pp. 103-126. \vspace*{0.1cm}

\item Lopes, S.R.C., B.P. Olbermann and V.A. Reisen (2004). ``A Comparison
of Estimation Methods in Non-stationary Arfima Processes". {\it
Journal of Statistical Computation and Simulation\/}, {\bf 74}(5),
pp. 339-347. \vspace*{0.1cm}

\item Lopes, S.R.C. and A. Pinheiro (2007). ``Wavelets for
Estimating the Fractional Parameter in Non-stationary ARFIMA
Process". In revision. \vspace*{0.1cm}

\item Lopes, S.R.C. (2007). ``Topics on Long-range Dependence".
In revision. \vspace*{0.1cm}

\item Maes, C., F. Redig, F. Takens, A.V. Moffaert and E. Verbitski
(1999). ``Intermittency and weak Gibbs states''. Institut voor
Theoretisch Fysica K.U. Leuven, B\'elgica. \vspace*{0.1cm}

\item Mandelbrot, B.B. (1997). {\it Fractals and Scaling in
Finance: Discontinuity, Concentration, Risk\/}. New York:
Springer-Verlag. \vspace*{0.1cm}

\item Olbermann, B.P. (2002). {\it Estima\c c\~ao em Classes de
Processos Estoc\'asticos com Decaimento Hiperb\'olico da Fun\c c\~ao
de Autocorrela\c c\~ao\/}. Ph.D. Thesis in the Mathematics Graduate
Program. Federal University of Rio Grande do Sul, Porto Alegre. URL
Address: www.mat.ufrgs.br/$\sim$slopes. \vspace*{0.1cm}

\item Olbermann, B.P., S.R.C. Lopes and V.A. Reisen (2006).
``Invariance of the First Difference in ARFIMA Models".  {\it
Computational Statistics\/}, {\bf 21}(3), 445-461. \vspace*{0.1cm}

\item Peng, C.K., S. Havlin, J.M. Hausdorff, J. Mietus, H.E.
Stanley and A.L. Goldberger (1996). ``Fractal Mechanisms and Heart
Rate Dynamics: Long-range Correlations and their Breakdown with
Disease". {\it Journal of Electrocardiology}, {\bf 28}, pp. 59-65.
\vspace*{0.1cm}

\item Peng, C.K., S.V. Buldyrev, A.L. Goldberger, S. Havlin, F.
Sciortino, M. Simons and H.E. Stanley (1992). ``Long-range
correlations in nucleotide sequences". {\it Nature\/}, {\bf 356},
pp. 168-170. \vspace*{0.1cm}

\item Percival, D.B. e A.T. Walden  (1993). {\it Spectral Analysis
for Physical Applications: Multitaper and Conventional Univariate
Techniques.} Cambridge: Cambridge University Press. \vspace*{0.1cm}

\item Pinheiro, A.S. and S.R.C. Lopes (2007). ``Bias Corrected
Wavelet Estimation for Fractional Parameter". In preparation.
\vspace*{0.1cm}

\item Reisen, V.A., B. Abraham and S.R.C. Lopes (2001) ``Estimation
of Parameters in the ARFIMA Processes: A Simulation Study". {\it
Communications in Statistics: Simulation and Computation}, {\bf
30}(4), pp. 787-803. \vspace*{0.1cm}

\item Reisen, V.A. and S.R.C. Lopes (1999). ``Some Simulations and
Applications of Forecasting Long-Memory Time Series Models''. {\it
Journal of Statistical Planning and Inference}, {\bf 80}(2), pp.
269-287.

\item Schuster, H.G. (1984). {\it Deterministic Chaos - An
Introduction}. Weinheim: Physik-Verlag.

\item Thaler, M. (1980). ``Estimates of the invariant densities of
endomorphism with indifferent fixed points''. {\it Israel Journal of
Mathematics}, {\bf 37}(4), pp. 303-313. \vspace*{0.1cm}

\item Wang, X.-J. (1989). ``Statistical Physics of Temporal
Intermettency". {\it Physics Review A\/}, {\bf 40}(11), pp.
6647-6661. \vspace*{0.1cm}

\item Young, L.-S. (1999). ``Recurrence times and rates of
mixing''. {\it Israel Journal of Mathematics}, {\bf 110}, pp.
153-188. \vspace*{0.1cm}

\end{description}

}

\end{document}